\documentclass[12pt,leqno]{amsart}
\usepackage[dvips]{graphics}
\usepackage{amssymb}

\numberwithin{equation}{section}

\newtheorem{thm}{Theorem}[section]

\newtheorem{lem}[thm]{Lemma}
\newtheorem{cor}[thm]{Corollary}
\newtheorem{prop}[thm]{Proposition}

\newtheorem{claim}[thm]{Claim}

\theoremstyle{definition}
\newtheorem{defn}{Definition}[section]
\theoremstyle{remark}
\newtheorem{rem}{Remark}[section]
\newtheorem{ex}[rem]{Example}

\newcommand{\tref}[1]{Theorem~\ref{#1}}
\newcommand{\cref}[1]{Corollary~\ref{#1}}
\newcommand{\pref}[1]{Proposition~\ref{#1}}

\newcommand{\lref}[1]{Lemma~\ref{#1}}

\newcommand{\Iso}{\mathrm{Iso}}

\newcommand{\R}{\mathbb{R}}
\newcommand{\Id}{\mathrm{Id}}

\begin{document}
\pagebreak


\title{THE DE RHAM DECOMPOSITION THEOREM FOR METRIC SPACES}

\author{Thomas Foertsch and Alexander Lytchak}
\address{Mathematisches Institut\\ Universit\"at Bonn\\
Beringstr. 1, 53115 Bonn, Germany\\}
\email{foertsch\@@math.uni-bonn.de \\ lytchak\@@math.uni-bonn.de}



\begin{abstract}
We generalize the classical de Rham decomposition theorem for Riemannian manifolds
to the setting of geodesic metric spaces of finite dimension.
\end{abstract}

\maketitle
\renewcommand{\theequation}{\arabic{section}.\arabic{equation}}
\pagenumbering{arabic}



\section{Introduction}
The direct product of metric spaces $Y$ and $Z$ is the Cartesian product 
 $X=Y\times Z$ with the metric given by
$d((y,z), (\bar y ,\bar z)) =\sqrt { d ^2(y,\bar y) + d^2 (z, \bar z)}$. Call a metric space
$X$ irreducible if for each decomposition $X=Y\times Z$ one of the factors $Y$ or $Z$ must
be a point.  It is a very  natural question if a given metric space has a ``unique'' decomposition
as a product of finitely many irreducible spaces.

 In general, no finite decomposition as a product of irreducible spaces may exist,
as an infinite product (for instance, a Hilbert cube) shows. On the other hand, there
is no uniqueness in general, even for some subsets of the Euclidean plane (\cite{four}, \cite{intui}).

 In the realm of Riemannian geometry this question is answered by the classical Theorem of de Rham 
(\cite{rham}). It says that for a simply connected, complete Riemannian manifold $M$ and each point 
$x\in M$,  subspaces of the tangent space $T_x M$ that are invariant under the action
of the holonomy group $Hol _x$  are in one-to-one correspondence  with the factors of $M$. As a consequence
one derives that each simply
connected complete Riemannian manifold $M$ admits a unique decomposition as a direct product
 $M= M_0 \times M_1 \times ...\times M_k$, where $M_0$ is a Euclidean space $\R ^m$ (possibly a point)
and all $M_i$,   $i=1,...,k$, are   irreducible Riemannian manifolds (with more than one point)
not isometric to the real line.   The uniqueness states that the factors are determined not only up to 
an abstract isometry, but that the $M_i$-fibers through a given point $x$ (i.e. the fiber 
$(M_i)_x :=P_i ^{-1} (P_i (x))$
 of the projection
$P_i  :M\to M_0 \times ...\times \hat M_i \times .... \times M_k$) are uniquely determined 
up to a permutation of
indices. 
Observe that the Euclidean space plays a special role, 
since it has many different decompositions as a product of real lines.

In \cite{eh} the statement about the uniqueness of the decomposition of $M$ was generalized to non simply connected
complete Riemannian manifolds by studying the action of the fundamental group of $M$ on the product 
decomposition of the  universal covering. 

Our main result presented in this paper is  a broad generalization of
de Rham's decomposition theorem. 
In order to state it precisely, recall that a geodesic in a metric space 
$X$  is an isometric embedding 
of an interval into $X$. A metric space is called geodesic if each pair of its points
is connected by a geodesic. 
We say that a metric space is affine if it is isometric to a linearly convex subset of a normed vector 
space. Given a metric space $X$ we define its affine rank, $\operatorname{rank}_{\operatorname{aff}}(X)$,
as the supremum over all topological dimensions of affine spaces that admit an isometric embedding
into $X$. Note that $\operatorname{rank}_{\operatorname{aff}}(X)$  is bounded above by the topological
dimension of $X$.

With this terminology our main result reads as follows.

\begin{thm} \label{main}
Let $X$ be a geodesic metric space of finite affine rank.
Then $X$ admits a unique  decomposition as a direct product
\begin{displaymath}
X=Y_0 \times Y_1 \times Y_2 \times ... \times Y_n ,
\end{displaymath} 
where $Y_0$ is a Euclidean space (possibly a point), and where the $Y_i$, $i=1,...,n$, 
are irreducible metric spaces not isometric to the real line nor to a point.
Thus, if there is another direct product decomposition $X=Z_0 \times Z_1 \times ...\times Z_m$ of
this kind then we have $m=n$ 
and there exists a permutation $s$ of $\{0, 1,...,n\}$
such that for each point $x\in X$ the $Y_i$-fiber through
$x$ coincides with the $Z_{s(i)}$-fiber through $x$ for all $i=1,...,n$. 
\end{thm}

\begin{rem}
Note that the Theorem \ref{main} cannot be a consequence of some kind 
of ``general nonsense'', as the examples of \cite{four} and \cite{intui}
demonstrate. Moreover, there is no similar decomposition
theorem in many  other categories. For instance, there are  finitely generated groups that have many completely
different decompositions as a direct sum   (cf.  \cite{baum}) and there are manifolds (for instance the Euclidean
space) that have  completely different decompositions as products of irreducible manifolds. 
\end{rem}

\begin{rem}
For compact subsets $X$ of a Euclidean space $\R^n$ the uniqueness of the decomposition
of $X$ was proved in \cite{mos}. 
\end{rem}

In the formulation of \tref{main} the Euclidean spaces again play a special role. As a particular case 
of \tref{main} 
(that is also an important step in the proof)  this special role can be expressed as
a funny rigidity statement:
 \begin{cor} \label{strike}
Let $X$ be a geodesic metric space of finite affine rank.  Assume that $X$ is decomposed in two
different ways as $Y\times \bar Y=X = Z\times  \bar Z$. Assume that the decompositions are 
transversal at some point 
$x\in X$, i.e. for the fibers $Y_x, \bar Y _x ,Z_x ,\bar Z _x$ we have:
$Y_x \cap Z_x = Y_x \cap \bar Z_x = \bar  Y_x \cap Z_x = \bar Y_x \cap \bar Z_ x = \{ x \}$. Then
$X$ is a Euclidean space $X=\R ^{2m}$.  
\end{cor}

Using \tref{main} we can analyze the group of isometries of a product space. To appreciate
the next result and the special role of Euclidean spaces, one should recall,
that the isometry group $\Iso (\R ^2)$of the plane $\R ^2$ is $3$-dimensional, and much larger than the 
two dimensional group $\Iso (\R ) \times \Iso (\R )$, i.e. the product $\R \times \R$ has many isometries
that do not respect the product structure. On the other hand we have:

\begin{cor} \label{isom}
Let $X$ be a geodesic space of finite affine rank and let $X=Y_0\times Y_1 \times ... \times Y_n$
be its product decomposition  as in \tref{main}. Denote by $\mathcal P$ the group
of all permutations $s\in \sigma _n$, such that $Y_i$ and $Y_{s(i)}$ are isometric for all  $i=1,...,n$.
 Then there is a natural exact sequence:
 \begin{displaymath}
1\to \Iso (Y_0)  \times \Iso (Y_1) \times ... \times \Iso (Y_n) \stackrel{i}{\to} \Iso (X) 
\stackrel{p}{\to} \mathcal P \to 1.
\end{displaymath} 
\end{cor}

\begin{rem}
 In the formulation of the results above and below, the equality sign in
$X=Y_0 \times ... \times Y_n$ or $X=Y\times \bar Y$  should be understood as a fixed isometry 
$I:X\to  Y_0 \times ... \times Y_n$ or $I:X\to Y\times \bar Y$. 
This fixed isometry then defines projections 
$P_i:X\to Y_i$ and $P^Y :X\to Y$, respectively. It also defines  $Y$-fibers, as
$Y_x =(P ^{\bar Y} )^{-1}(P^{\bar Y} (x))$. 
\end{rem}

We are going to explain the basic idea behind the proof of our main theorem. 
The idea is to find
a class $\mathcal C$ of metric spaces that is small enough, such that one can prove the main theorem
for this class by some direct means, and that is large enough such that each metric space $X$
as in \tref{main} can be approximated by elements of $\mathcal C$ in some suitable sense.

We use the class $\mathcal C$ of affine metric spaces of finite dimension. To ``approximate''
general spaces by elements of $\mathcal C$ we prove that for each metric space $X$, 
each maximal affine subset $C$ of $X$ is ``rectangular'', i.e. for each decomposition
$X=Y\times \bar Y$ we have $C=P^Y (C) \times P^{\bar Y} (C)$. This claim is a direct consequence
of the following technical result:

\begin{prop} \label{invar}
Let $X=Y\times \bar Y$ be a direct product. If a subset $C\subset X$ is affine, then so is
the projection  $P^{Y} (C)\subset Y$. 
\end{prop} 

\begin{rem}
 Essentially, the statement of \pref{invar} is that for $C\subset X$ as above,
parallel  linear geodesics in $C$ have the same slope with respect to the product
decomposition of $X$.  If $C$ is a normed vector space (i.e. if each geodesic in $C$ is part
of an infinite geodesic) then this claim is easy to verify and \pref{invar} for the case of normed  
vector spaces $C$ was already proven in \cite{fs}. In our more general case, the result is more subtle
and relies on the infinitesimal considerations in \cite{hitz}. 
 In \cite{fs} it is shown by an example, that
\pref{invar} becomes wrong if one replaces the word ``affine'' by ``convex subset of a Euclidean space''.
\end{rem}

 Now we explain the heart of the proof of \cref{strike}. Under the assumptions of \cref{strike}
we choose a largest affine subset $C$ of $X$ that contains $x$.  We know that $C$ is rectangular,
hence it has two mutually transversal decompositions, as well. Assuming that we already
know the result for affine spaces, we deduce that $C$ is a Euclidean space $C=\R ^{2m}$. It
remains to prove that $C=X$.  First we observe that each rectangular subset of $X$
that contains $C \cap Y_x$ must contain the whole subset $C$ (by easy linear algebra). 
Now take an arbitrary point 
$z\in \bar Y_x$. Then $z$ and $C\cap Y_x$ ``span'' a flat subset of $X$. Let $C_0$ be a largest
affine subset of $X$ that contains this flat one. Then $C_0$ is rectangular, hence
it contains $C$ and by maximality of $C$ we have $C=C_0$. Thus $Y_x \subset C$.
 Interchanging the roles of $Y$ and  $\bar Y$ we deduce the result.

{\bf Outline of the paper:} In Section \ref{sec-prelim} we discuss basic facts about metric products.
In Section \ref{sec-affine} we discuss affine metric spaces and prove \pref{invar}. In
Section \ref{sec-cons} we draw some direct consequences of the equality of slopes in an affine subset
of a product.  
In Section \ref{sec-inter} and Section \ref{sec-lin} we prove that for different decompositions
$Y\times \bar Y =X=Z\times \bar Z$ and each point $x\in X$ the intersection $Y_x \cap Z_x$ is a factor
of $Z_x=Z$.  In Section \ref{sec-reduce} we use this observation to reduce \tref{main} to \cref{strike}. 
Finally, in Section \ref{sec-ban} we prove \cref{strike}.

\section{Preliminaries}
\label{sec-prelim}


\subsection{Notations and basic observations}
\label{subsec-notations}
By $d$ we will denote distances in metric spaces without an extra reference 
to the space. 

 For a direct product $X=Y\times \bar Y$ we will use the following notations.

By $P^Y:X\to Y$ and $P^{\bar Y} :X\to \bar Y$ we will denote the canonical 
projections to the factors. (In fact the equality $X=Y\times \bar Y$ just
means that two maps $P^Y:X\to Y$ and $P^{\bar Y} \to \bar Y$ are given
such that $d^2(x,z)=d^2 (P^Y(x),P^Y(z)) +d^2 (P^{Y} (x),P^{\bar Y} (z))$
holds for all $x,z\in X$).

  For a point $x\in X$ we call the subset $(P^{\bar Y})^{-1} (P^{\bar Y} (x))$
the $Y$-fiber through $x$ and denote it by $Y_x$. The restriction
$P^Y :Y_x \to Y$ is an isometry and we will sometimes identify $Y_x$ with
$Y$ via this isometry. The composition $P^{Y_x}$
of $P^Y:X\to Y$ and the inverse of $P^Y:Y_x\to Y$ is the natural projection
of $X$ to $Y_x$. It sends a point $z\in X$ to the unique intersection
of $\bar Y_z$ and $Y_x$.  Each restriction $P^{Y_x} :Y_z \to Y_x$ is an isometry.
 For all $x,z \in X$ we have
$d^2 (x,z) =d^2 (x, P^{Y_x} (z))+ d^2(P^{Y_x} (z),z)$.


\subsection{Recognition of products} \label{recogn}
 Assume on the other hand that a space $X$ is given as a 
 union $X=\cup _{i\in J}  Y_i$ (We do not assume the union to be disjoint, 
moreover, $Y_i$ and $Y_j$ may coincide for different $i$ and $j$). Assume that
for all $i,j \in J$ a map $P_{ij} :Y_i \to Y_j$ is given such that
for all $i,j,k \in J$ we have $P_{ij} \circ P_{ji} =Id$ and
$P_{jk} \circ P_{ij} =P_{ik}$. Furthermore, assume that 
for all $i,j \in J$, all $x\in Y_i$ and $\bar x \in Y_j$ we have
$d^2(x,\bar x) = d^2 (x,P_{ij} (x)) + d^2 (P_{ij} (x), \bar x)$.

 We are going to show that $X$ splits as a product with fibers $Y_i$.

For all $x,\bar x$ as above we have:
\begin{displaymath}
d^2 (x,\bar x) = d^2 (x,P_{ij} (x)) + d^2 (P_{ij} (x), \bar x) =
d^2(\bar x, P_{ji} (\bar x) ) + d^2 (P_{ji} (\bar x), x)
\end{displaymath} 
and 
\begin{displaymath}
d^2 (P_{ij} (x), P_{ji} (\bar x)) 
= d^2 (P_{ij} (x), \bar x )+
d^2 (\bar x, P_{ji} (\bar x)) = d^2 (P_{ij} (x),x) + d^2 (x, P_{ji} (\bar x)). 
\end{displaymath} 

Subtracting these equalities from another we obtain
$d^2 (x,P_{ij} (x)) =d^2 (\bar x, P_{ji} (\bar x))$.
Therefore  $d(x,P_{ij} (x)) =d(\bar x,P_{ji} (\bar x))$.
In particular, $P_{ij}$ is an isometry and we have
$d(Y_i,Y_j) =d(x,P_{ij} (x))$ for each $x\in Y_i$.

Therefore, the $Y_i$ define a so called equidistant
decomposition of $X$ and $d(i,j):= d(Y_i,Y_j)$ defines a pseudo metric on
$J$, where two indices have distance $0$ if and only if they define
the same subset of $X$. We identify equal fibers and may assume that different
fibers are disjoint, i.e. that $J$ is a metric space.

For all $i,j \in J$ and all $x\in Y_i,\bar x \in Y_j$ we have:
$d^2 (x,\bar x) = d^2 (P_{ij} (x), \bar x) +d^2 (i,j)$.

Fix now a fiber $Y_o$ for some $o\in J$ and consider the map
$P:Y_0\times J \to X$ given by $P(y,i) =P_{oi} (y)$.
This map is surjective, by assumption. We claim that it is an isometry.

Indeed, for all $y,\bar y \in Y_o$ and all $i,j \in J$ we have 
$d^2 (P_{oi} (y), P_{oj} (\bar y)) = d^2 (i,j) + 
d^2 (P_{ij} (P_{oi} (y)),P_{oj} (\bar y))$.

But $P_{ij} (P_{oi} (y)) =P_{oj} (y)$ and
$d(P_{oj} (y),P_{oj} (\bar y)) =d(y,\bar y)$. 
This finishes the proof.

\subsection{Intersections of different fibers}
Let $Y\times \bar Y=X=Z\times \bar Z$ be two decompositions of a space
$X$. Fix a point $x\in X$ and set $F_x =Y_x \cap Z_x$. The following
lemma together with the preceding subsection suggests that
$F_x$ has good chances to be a factor of $Z_x$.

\begin{lem} \label{interbase}
For each point $p\in \bar Y_x$ and each point $q \in F_x$ we 
have $d^2 (P^Z (p),P^Z(q))=
d^2 (P^Z (p),P^Z (x)) + d^2 (P^Z(x),P^Z(q))$.
\end{lem}

\begin{proof}
Identify $Z$ with $Z_x$ and denote by $\tilde p$ the projection of
$p$ onto $Z_x$. Observe that $x$ and $q$ are already in $Z_x$.

We have $d^2 (q,p)=d^2 (q,\tilde p) +d^2 (\tilde p,p)$ and
$d^2 (x,p)=d^2 (x,\tilde p)+ d^2 (\tilde p,p)$. Since $q\in Y_x$ and
$p\in \bar Y_x$  we have $d^2 (q,p)=d^2 (q,x) + d^2 (x,p)$.
 We  insert the second equality in the third and the third in the first one
and get $d^2 (\tilde p,q) = d^2 (\tilde p,x) + d^2 (x,q)$.
\end{proof}

\subsection{Geodesics in products} Recall that geodesics (if not otherwise stated)
are parameterized by the arclength.
A subset $C$ of a geodesic metric space $X$
is called convex (totally convex, resp.) if for each pair of points
$y,z\in C$ there is some geodesic in $C$ between these two points
(if each geodesic between $y$ and $z$ is contained in $C$, resp.).

Each geodesic $\gamma$ in a product $X=Y\times \bar Y$ has the form
$\gamma (t) = (\eta (at),\bar \eta   (\bar a t))$ for some geodesics
$\eta$ in $Y$ and $\bar \eta $ in $\bar Y$ and some real numbers $a,\bar a$
with $a^2 +\bar a ^2 =1$.  The numbers $a$ and $\bar a$ are called the slopes
of the geodesic $\gamma$ with respect to $Y$  and to $\bar Y$, respectively.
Note that   the slope of $\gamma$ with respect to $Y$ is $0$ if and only if
$\gamma$ is contained in some $\bar Y$-fiber.

A product $X=Y\times \bar Y$ is a geodesic space if and only if the factors
$Y$ and $\bar Y$ are geodesic. In this case each fiber $Y_x$ is totally convex
in $X$. The fact that geodesics project to geodesics implies that for each
convex subset $C$ in a product $X=Y\times \bar Y$ the projection 
$P^Y (C)$ is a convex subset  of $Y$.

\subsection{Groups of isometries} We are  going to deduce \cref{isom} from \tref{main}.
Thus let $X=Y_0\times Y_1 \times ... \times Y_n$ be as in \cref{isom} and let
$g:X\to X$ be an isometry.  By \tref{main} the isometry $g$ must induce a permutation
$s:\{ 0,1,...,n\} \to \{ 0,1,...,n \}$ such that $g ((Y_i)_ x) =(Y_{s(i)} )_{g(x)}$.
The restriction $g:(Y_i)_x \to  (Y_{s(i)} )_{g(x)}$ must be an isometry, hence $Y_i$
and $Y_{s(i)}$ must be isometric. In particular, $s(0)=0$. The assignment $g\to s$ 
is a homomorphism $p:\Iso (X) \to \mathcal P$.  Interchanging isometric factors by
some fixed isometry, one sees that $p$ is surjective.

 The map $i:\Iso (Y_0) \times ... \times \Iso (Y_n) \to \Iso  (Y_0\times ...\times Y_n)$
is given by $i (g_0,g_1,...,g_n) (y_0,y_1,...,y_n)=(g_0(y_0),g_1(y_1),...,g_n(y_n))$.
It is a well defined injection and satisfies $p\circ i =1$. Let now $g$ be an element
in the kernel of $p$. Fix a point $x\in X$. Then $g$ induces isometries $g:(Y_i)_x \to
(Y_i)_{g(x)}$.  Identifying $(Y_i)_x$ and $(Y_i)_{g(x)}$ with $Y_i$ via the projection
$P^{Y_i}$ we obtain an isometry $g_i :Y_i \to Y_i$. Set $\tilde g = i (g_0,g_1,...,g_n)$.
We claim $g= \tilde g$. Consider the isometry $h =g \circ \tilde g ^{-1}$.
We have $h(x)=x$ and the restriction of $h$ to $(Y_i)_x$ is the identity for all $i$.
For each point $\bar x \in X$ we must have $P^{(Y_i)_x} (\bar x) = P^{(Y_i)_x} (h(\bar x))$
for all $i$. But this means $P^{Y_i} (\bar x)=  P^{Y_i} (h(\bar x))$, hence $\bar x= h(\bar x)$.
This shows $h=Id$ and proves that the sequence is exact.


\section{Affine spaces}
\label{sec-affine}
\subsection{Definitions} We are going to prove \pref{invar} in this section.
If we want to prove that projections of affine subsets are affine,
we first have to find linear geodesics in the projection. Thus we need to work
with a distinguished class of geodesics.

\begin{defn}
Let $X$ be a metric space. A bicombing $\Gamma$ on $X$ assigns to each pair of points
$(x,y)\in X\times X$ a geodesic $\gamma_{xy}$ connecting $x$ to $y$, such that
$\gamma_{xy}=\gamma_{yx}$ as sets (where $\gamma_{xy}$
and $\gamma_{yx}$ have opposite orientation as curves) and such that for each $m\in \gamma_{xy}$ 
we have $\gamma_{my}\subset \gamma_{xy}$.
\end{defn}

A space with a bicombing is per definitionem a geodesic metric space.
Given a space with a bicombing $\Gamma$ we will call the geodesics assigned
by $\Gamma$ the special geodesics.

\begin{ex}
Let $X$ be a geodesic metric space in which there is only
one geodesic  connecting each pair of points. Then $X$ has a (unique)
natural bicombing. 
\end{ex}

A map $I:X\to Y$  between two spaces with bicombings $\Gamma$ and $\Gamma'$ is called
affine if it sends special geodesics to special geodesics parameterized
proportionally to the arclength. Such a map  $I$ is called
a $(\Gamma ,\Gamma')$-isometric embedding ($(\Gamma ,\Gamma')$-isometry, resp.) if 
it is an isometric embedding (an isometry, resp.) between the underlying metric
spaces.

A subset $Y$ of a space $X$ with a bicombing $\Gamma$ will be called 
$\Gamma$-convex if for each pair of points $x,y\in Y$ the geodesic
$\gamma _{xy}$ is contained in $Y$.  Observe that the image of each
$(\Gamma ,\Gamma')$-isometric embedding is a $\Gamma'$-convex subset and that
each $\Gamma$-convex subset inherits a natural bicombing.

 If $(X_i,\Gamma _i)$ are spaces with bicombings for $i=1,2$, then
the direct product $X=X_1\times X_2$ has a unique bicombing such that
the projections $P^{X_i} :X\to X_i$ are affine.

 For us, the main example will be the following. 
Let $V$ be a normed vector space. Then the linear geodesics (i.e.
geodesics of the form $\gamma (t)=v +tw$) define a bicombing on $V$.
Normed vector spaces will always be equipped with this particular
bicombing.

 A subset $C$ of a normed vector space $V$ is $\Gamma$-convex if and only
if it is linearly convex (i.e. for each $x,y \in C$ and each $t\in [0,1]$
the point $tx + (1-t) y$ is in $C$).

 The notion of affine maps in this setting coincides with the usual one.
First, let $V,W$ be normed vector spaces. Then each linear map $A:V\to W$ 
is affine. On the other hand, let $F:V\to W$ be an affine map
with $F(0)=0$. Then $F(tx)=tf(x)$  and 
$F(x+y /2) =F(x+y)/2$, for  each $t\in \R$ and all $x,y\in V$, hence $F$
is linear in this case.   By adding a translation
we deduce that a map $F:V\to W$ is affine if and only if it has the form
$F(v)=w_0 +A(v)$ for some $w_0 \in W$ and some linear map $A:V\to W$.

 Let now $X$ be a linearly convex subset of a normed vector space
$V$ and let $f:X\to W$ be an affine map into a normed vector space
$W$.  We first assume that $0\in X$ and that $f(0)=0$. 
Denote by $C$ the cone over $X$, i.e.  $C=\{ tx | t\geq 0,x\in X\}$.
Then $C$ is linearly convex in $V$ and $\bar f  (tx)=t f(x)$ is a well
defined extension of $f$ to $C$ that is again affine. 
Denote by $V_0$ the linear hull of $X$. Then $V_0 =\{ x-y |x,y\in C \}$.
The map $\tilde f (x-y)=\bar f (x) -\bar f (y)$ is well defined and affine.
By the above observation the map $\tilde f$ is linear and can be extended
to a linear map $F:V\to W$, by linear algebra.  We deduce,
that a map $f:X\to W$ is affine if and only if it is the restriction
of an affine map $F:V\to W$.

\begin{defn}
 We call a space $X$ with a bicombing $\Gamma$ affine if  
there is a $\Gamma$-isometric embedding 
$I:X\to V$ into a normed vector space.
\end{defn}

Thus a metric space is affine in the sense of the introduction if and
only if it has a bicombing $\Gamma$ such that
$(X,\Gamma )$ is affine.

\begin{rem} 
Each metric space $X$ has an isometric embedding into a Banach space.
If $X$ is geodesic then the image of $X$ is a convex subset of $X$.
In order to make the above definition not trivial it is necessary 
to distinguish (as we did) between convex and linear convex subsets.   
\end{rem}


\subsection{Three points characterization} 

\begin{prop} \label{prop-three-points}
Let $X$ be a space with a bicombing $\Gamma$. If for all points $x,y,z \in X$
there is a $\Gamma$-convex subset $C_{x,y,z}$ of $X$ 
that is affine and contains
the three points $x,y$ and $z$, then $X$ is affine.  
\end{prop}

\begin{proof}
Fix a point $o$ in $X$.
 Consider the space $Y=X\times [0,\infty )$
 and identify  points $(x,t)$ and $(z,s)$ if we have 
$t\cdot  d(o,x) = s \cdot d(o,z)$ and the geodesics $\gamma _{ox}$ and
$\gamma _{oz}$ initially coincide. Moreover we identify the subsets
$X\times \{0 \}$ and $\{ o\} \times [0,\infty )$ with a unique point
in $Y$ that we  denote by $0$ and call the origin.

Denote by $Z$ the arising set. 
On $Z$ we define a metric by setting $d((x,t),(y,s)):= \frac 1 \epsilon  
d(\gamma _{ox} (\epsilon t), \gamma _{oy} (\epsilon s ))$ for a sufficiently small
number $\epsilon >0$. This quantity is well defined 
(i.e. does not depend on $\epsilon$) as one sees by making the computations
 in the linear subset $C_{o,x,y}$ as in the assumptions.
 Moreover $d$ is a metric as one derives from the triangle inequality in $X$.

We have a multiplicative operation of $[0,\infty )$ on $Z$ by
$\lambda \cdot (x,t) := (x, \lambda t )$.  For 
$z_1,z_2 \in Z $ we have by definition 
$d(\lambda z_1, \lambda z_2 ) = \lambda d(z_1,z_2)$.

Next the space $X$ has an isometric embedding $I$ into $Z$ defined
through $I(x) = (x,1)$. We identify $X$ with its image in $Z$.
By our definition we see that for each finite subset $D$ of $Z$ there is
some $\lambda>0$ such that $D$ is contained in $\lambda X$.

Taking two points $z_1,z_2\in Z$ we find some $\lambda >0$ and points
$x_1,x_2 \in X$ with $\lambda x_i =z_i$. We extend the bicombing from
$X$ to $Z$ by letting the geodesic $\gamma _{z_1z_2}$ be the (reparameterized)
curve $\lambda \gamma _{x_1x_2}$. Again considering the triangle 
$C_{x_1,x_2,o}$ we see that this definition does not depend on the choice
of $\lambda $ (i.e. on the choice of $x_1$ and $x_2$) and that this 
is indeed a bicombing on $Z$.  Moreover, for three points 
$z_1,z_2,z_3\in Z$ we define a $\Gamma$-convex set $C_{z_1,z_2,z_3}$ by
choosing $x_i$, $i=1,2,3$, and $\lambda $ as above, setting
$C_{z_1,z_2,z_3} := \lambda C_{x_1,x_2,x_3}$. It follows that $C_{z_1,z_2,z_3}$ is affine.

Hence $Z$ again has the same $3$-point property as $X$ and,
since $X$ isometrically embeds in $Z$,
it is enough to prove that $Z$ is affine.

For this purpose we define an addition on $Z$ in the following way. For $x,y \in Z$ we set
$x+y := 2m$ where $m$ is the midpoint of $\gamma _{xy}$. By definition we
have $x+0=x$ and $x+y=y+x$ for all $x,y\in Z$.
Considering the subset $C_{o,x,y}$ we see
that for all $x,y\in Z$ we have $\lambda (x+y) =\lambda x +\lambda y$ for each $0\leq \lambda \leq 1$
and therefore for each $\lambda \geq 0$.

Let $x,y\in Z$ be arbitrary. Considering the affine space
$C_{o,2x,2y}$ we see that for all $0\leq t \leq 1$ we have
$tx+ (1-t)y = \gamma _{xy} (t d(x,y))$.

Take now three points $x,y,z \in Z$. Then  
$u:=\frac 1 3 (((x+y) +z)= \frac  2 3 (\frac 1 2 (x+y) ) + \frac 1 3 z$. 
This shows that $u \in C_{x,y,z}$. The same is true for
$\bar u:= \frac 1 3 (x+ (y+z))$. Making the computations in
$C_{x,y,z}$ we derive that $u=\bar u$. This implies that the addition
is associative.

 Finally considering for three arbitrary points $x,y,z \in Z$ the 
subset $C_{x,y,z}$ we see that $\frac {x+y} 2 = \frac {x+z} 2$ implies 
$y=z$.

Hence $Z$ with the defined addition is a commutative semi-group 
in which $x+y=x+z$ implies $y=z$. Let $W$ denote the abelian group defined by
$Z$, i.e. $W$ is defined as the set of all pairs $(x,y)\in Z^2$ modulo the 
equivalence relation $(x,y)\sim (\bar x, \bar y)$ if and only if  $x+\bar{y}=\bar x + y$.
The addition on $W$ is defined by $(x,y)+(\bar x, \bar y)=(x+\bar x, y+\bar y)$. 

The set $W$ with this addition is an abelian group and $Z$ has the canonical
embedding $i:Z\to W$ defined by $i(x):=(x,0)$. The operation of $\R ^+$ on 
$Z$ extends to an operation on $W$ by setting 
$\lambda (x,y) :=(\lambda x,\lambda y)$. Finally we set 
$(-1) \cdot (x,y) := (y,x)$ and for each positive $\lambda$ and each $w\in W$
we set $(-\lambda) \cdot w:= (-1) \lambda w$.
 Taking this together we see that $W$ is a vector space over $\R$.
Moreover, the embedding $i:Z\to W$ sends $\Gamma$-geodesics to linear intervals.

 Finally we define the function $|\cdot |$ on $W$ by setting $|(x,y)|:=d(x,y)$.
It remains to prove that $|\cdot|$ is well defined and that it is a norm.

 To see that it is well defined let $x,y,\bar x, \bar y \in Z$ be such that
$(x,y)\sim (\bar{x},\bar{y})$, i.e.
$x+\bar y =\bar x + y$.  Then the midpoint $m$ of $\gamma _{x\bar y}$
is also the midpoint of $\gamma_{\bar x y}$. Take the midpoint $p$ between
$x$ and $y$. Considering the space $C_{x,y,\bar x}$ we see that 
$\gamma _{p\bar x}$ intersects $\gamma _{x\bar y}$ in some point $p_1$ that is 
between    $x$ and $m$. Similarly $\gamma _{p\bar y}$ intersects $\gamma_{ym}$
in some point $p_2$. Therefore the subspace $C_{p,\bar x,\bar y}$ contains 
$p_1,p_2$ and $m$.  Hence for some small $\epsilon >0$ we can find 
points $q_1$ and $\bar q_1 $ on $\gamma _{p_2 m}$ and $\gamma _{m\bar x}$ 
and points $q_2$ and $\bar q_2$ on $\gamma _{p_1m}$ and $\gamma _{m \bar y}$
such that $d(m,q_1)=d(m,\bar q_1)= \epsilon d(m,y)$ and 
$d(m,q_2)=d(m,\bar q_2) =\epsilon d(m,x)$. 

In $C_{p,\bar x,\bar y}$ we deduce from this that 
$d(q_1,q_2)=d(\bar q_1,\bar q_2)$. Now looking at the space
$C_{x,y,m}$ ($C_{\bar x, \bar y,m}$, resp.) we see that 
$d(q_1,q_2) =\epsilon d(x,y)$ 
($d(\bar q_1,\bar q_2) =\epsilon d(\bar x, \bar y)$, resp.). This shows that 
$|(x,y)|=|(\bar x,\bar y)|$ and $|\cdot|$ is well defined. 

From the definition we immediately conclude that $|\lambda w| =|\lambda| |w|$ for
each $\lambda \in \R$ and each $w\in W$.  Finally, given two points 
$w=[(x,y)]$ and $v=[(\bar x,\bar y)]$ in $W$ we can write $w=[(x+\bar x,y+\bar x)]$
and $v=[(\bar x +x,\bar y +x)]$, i.e. we may assume that $x=\bar x$.
In this case one deduces from $C_{x,y,\bar y}$ that 
$d(x,y)+d(x,\bar y) \geq 2 d(x,\frac {y+ \bar y} 2)$ which implies
$|v|+ |w| \geq |v+w|$. 

This finishes the proof.
\end{proof}

 As a corollary we obtain:

\begin{cor} \label{maxexist}
Let $X$ be a metric space and let $Y$ be an affine subset of $X$
with a bicombing $\Gamma$. Then there is a maximal affine subset 
$Y'$ of $X$ that contains $Y$ and the bicombing of which extends $\Gamma$.
\end{cor}

\begin{proof}
By Zorn's lemma it is enough to prove that for a chain $Y_i$ of affine subsets
such that the bicombing of $Y_i$ extends the bicombing of $Y_j$
for $i\geq j$ their union $Y'$ is affine.

 This union $Y'$ has  a natural bicombing that extends the bicombings of all
$Y_i$. The three points property from \pref{prop-three-points} is satisfied by $Y'$, since
it is satisfied by all $Y_i$. Thus $Y'$ is affine.
\end{proof}


\subsection{Invariance under products}
Since a projection of a product onto a factor sends geodesics
to geodesics parameterized proportionally to the arclength,
\pref{invar} is a direct consequence of the following result:

\begin{prop} \label{vary1}
Let $(X,\Gamma)$ be an affine metric space with a bicombing. Let
$f:X\to Y$ be a surjective continuous map that sends each special geodesic
to a geodesic parameterized proportionally to the arclength. Then $Y$
is affine. More precisely, there is a unique bicombing $\Gamma '$ on
$Y$, such that $(Y,\Gamma ')$ is affine and such that
$f:(X,\Gamma )\to (Y,\Gamma ')$ is affine. 
\end{prop}

  Before embarking on the proof we will cite an important special case.
Assume namely that $f$ is bijective. Then the images of special geodesics
(reparameterized to have speed $1$) define a (unique) bicombing 
$\Gamma '$ on $Y$, such
that $f$ is affine with respect to this bicombing. The inverse 
$f^{-1} : (Y,\Gamma ') \to (X,\Gamma )$ is an (a priori not continuous) 
 affine equivalence. In this
case Theorem 1.3 of \cite{hitz} says that $(Y,\Gamma ')$ is in fact affine.
Thus Proposition \ref{vary1} is in fact the extension of Theorem 1.3 of \cite{hitz} to
the non-injective case.

\begin{proof}
 First observe that the restriction $f:C\to f(C)$ of $f$ to
a $\Gamma$-convex subset of $X$ again satisfies the assumption.

 If a required bicombing $\Gamma '$ exists we must have
$\gamma _{f(x_1)f(x_2)} =f(\gamma _{x_1x_2})$ (as sets)
for all $x_1,x_2 \in X$. In particular, if such a bicombing exists,
it is unique.
 We only need to prove that the above assignment is well defined and that
$Y$ equipped with this bicombing is affine.

Assume first that $X$ is a triangle, i.e. $X$ is the linearly convex
hull of three points in  a two dimensional
normed vector space $V$. If the triangle is degenerate, i.e. a point
or an interval, then $Y=f(X)$ is a point or an interval as well
and the conclusion is clear. Thus we may assume that $X$ is non-degenerate.
 Denote by $X_0$ the set of inner points of $X$ (with respect to $V$).
Observe that $X_0$ is linearly convex and dense in $X$.
Set $Y_0=f(X_0)$.
 Each fiber of $f$ in $X_0$ is linearly convex. More precisely, it
is the intersection of $X_0$ with an affine subspace of $X$. Hence there are 
three cases.

1) There is a point $x\in X_0$ such that $f^{-1} (f(x))=x$.
Denote by $U$ the set of all such points and let $O$ be a connected
component of $U$. 

 For a point $z\in X_0$ and a unit vector  $v\in V$ we set $|v|_z$ to
be the speed of the geodesic $\gamma (t)= f(z+tv)$. By continuity
of $f$ the function $|v|_z$ is continuous in $z$ and in $v$. Moreover,
a point $z$ is in $U$ if and only if $|v|_z >0$ for all unit vectors $v\in V$. By continuity of
$|v|_z$ the subset $U$ (and therefore $O$) is open.

 On each convex subset $O_1$ of $O$ the restriction $f:O_1 \to f(O_1)$
is bijective. Thus one can apply \cite{hitz}, and 
 \pref{vary1} is true for this restriction.
 Hence $f:=O_1 \to f(O_1)$ is the restriction 
of an affine map between normed vector space in this case. Thus
$|v|_{z_1} =|v|_{z_2}$ for all unit vectors $v\in V$ and all $z_1,z_2 \in O_1$.
By connectedness and local convexity of $O$ we deduce that 
$|v|_z$ is constant on $O$. Hence, by continuity of $|v|_z$, at each
boundary point $z$ of $O$ we still have $|v|_z >0$ for all unit vectors $v\in V$. 
By connectedness of $X_0$ this implies $X_0=U$. Hence $F:X_0 \to Y_0$ is
bi-Lipschitz in this case. Hence the continuous extension $F:X\to F(X)$ is bi-Lipschitz
as well.
Thus we can apply \cite{hitz} again and obtain the validity of the statement in this case.

2) Each fiber of $f$ is an interval. In this case one can rearrange
each triple of fibers $I_1,I_2,I_3$, such that $I_2$ is between
$I_1$ and $I_2$, i.e. for each point $x_1 \in I_1$ and $x_3 \in I_3$
the geodesic between $x_1$ and $x_3$ intersects $I_2$. 
This implies that for each triple of points in $Y_0$ one can rearrange them
such that one point is on a geodesic between the other two. But this implies
that $Y_0$ is an interval (i.e. a subset of a real line). 
In this case $Y$ is an interval, too (since $Y_0$ is dense in $Y$)
and we are done.

3) There is only one fiber of $f$. In this case $Y_0=f(X_0)$ is a point.
Hence $Y$ is  a point and there is nothing to prove.

 Therefore the statement is true if $X$ is a triangle. 

Let now $X$ be arbitrary.
 To prove that the bicombing
is well defined choose points $x_i,z_i \in X$, for $i=1,2$ with
$f(x_i)=f(z_i)$. Considering the triangles spanned by $x_1,x_2,z_1$
and $x_1,z_1,z_2$ and using the fact that the statement is true for triangles
we deduce that $f(\gamma _{x_1x_2})=f(\gamma _{x_1z_2} )=f(\gamma _{z_1z_2})$.
This shows that the bicombing on $Y$ making $f$ an affine map is well 
defined.

 For an arbitrary triple $y_1,y_2,y_3 \in Y$ choose $x_1,x_2,x_3 \in X$
with $f(x_i)=y_i$. Let $C$ be the triangle spanned by $x_i$. Then 
$f(C)$ is an affine $\Gamma '$-convex subset that contains the points $y_i$.
Applying \pref{prop-three-points} we deduce that $Y$ is in fact affine. 
\end{proof}

\section{First applications}
\label{sec-cons}

\subsection{Rectangular subsets} 
We will call a subset $S$ of a metric space $X$ rectangular with respect to the
product decomposition $X=Y\times \bar Y$ of $X$ if $S=P^{Y} (S) \times P^{\bar Y} (S)$ holds.
We will say that $S$ is a rectangular subset of $X$ if it is rectangular with respect
to each product decomposition of $X$.

\begin{ex}
If $X$ is irreducible, then each subset $S$ of $X$ is rectangular. If $X=\R^n$, then
the only rectangular subsets of $X$ are the whole space, the empty set and subsets with
only one point.
\end{ex}

Let now  $X$ be a metric space and let $C$ be a maximal affine
subset of $X$ (i.e.  $C$ with some bicombing $\Gamma$ is affine and 
there is no larger affine subset of $X$ the bicombing of which extends that of
$C$).  If $X$ is    decomposed as a direct product 
$X=Y \times \bar Y$, then due to \pref{invar} the projections
 $P^Y(C)$ and $P^{\bar Y} (C)$ are affine and therefore
so is the product $\tilde C=  P^Y(C) \times P^{\bar Y} (C) \subset
Y \times \bar Y$. Moreover,  the restrictions of the projections
to $C$ are affine (\pref{vary1}), thus the natural bicombing 
 of $ P^Y(C) \times P^{\bar Y} (C)$ extends the bicombing of $C$.
By maximality of $C$ we deduce $C= P^Y(C) \times P^{\bar Y} (C)$.
This shows:

\begin{cor} \label{maxrect}
Let $X$ be a metric space and $C$ be a maximal affine subset of $X$.
Then $C$ is a rectangular subset of $X$.
\end{cor}

 Since the dimension of a product of two affine spaces equals the 
sum of the dimension of the factors we get (compare \cite{fs}):

\begin{cor} \label{additive}
 For each decomposition $X=Y\times \bar Y$ we have 
\begin{displaymath}
\operatorname{rank}_{\operatorname{aff}} (X) = \operatorname{rank}_{\operatorname{aff}} (Y) + 
\operatorname{rank}_{\operatorname{aff}} (\bar Y).
\end{displaymath}
\end{cor}


\subsection{Equality of slopes} \label{equalslope}
From \pref{vary1} we deduce that if 
$C$ is an affine subset of a product
$X=Y\times \bar Y$ then each pair of (linear) geodesics in $C$ that are  
parallel in $C$ have the same slope with respect to $Y$ (and to $\bar Y$).

From this we derive the following.
Let $X=Y\times \bar Y$ and $X=Z\times \bar Z$ be two decompositions
of a geodesic metric space $X$. Let $\gamma$ be a geodesic in $Y$.
Then the slope of the geodesic $\gamma \times \{ \bar y \} \subset X$
with respect to $\bar Z$ does not depend on the point $\bar y \in \bar Y$.
Indeed, for different points $\bar y_1, \bar y_2 \in \bar Y$ choose
a geodesic $\eta$ connecting them. Then $\gamma \times \eta \subset X $
is a flat rectangle and the $\gamma \times {\bar y _i}$ are
parallel sides of it.  

 This observation allows us to speak of the slope of a geodesic
$\gamma \in Y$ with respect to $\bar Z$. 
Observe that this slope is $0$ if and only if for the endpoints $y_1,y_2$
of $\gamma$ and some (and therefore each) point $\bar y \in \bar Y$ 
the points $(y_i,\bar y) \in X$ are contained in the same $Z$-fiber.

 In particular, we see that if $Y_x$ is contained in $Z_x$ for some point
$x\in X$ then $Y_{\bar x}$ is contained in $Z_{\bar x}$ for each point 
$\bar x \in X$.


\section{Intersections and Projections}
\label{sec-inter}

Let $X$ be a geodesic metric space with two decompositions
$Y\times \bar Y =X =Z\times \bar Z$. Let $x\in X$ be a point and set
$F_x = Y_x \cap Z_x$. Consider $T=T_x := P^Y (F_x ) \times \bar Y \subset
Y\times \bar Y =X$.  Then we have:

\begin{lem} \label{lem-proj}
In the notations above the image $P^Z (T)\subset Z$ splits as
$P^Z(T)= P^Z(F_x) \times P^Z (\bar Y _x)$.
\end{lem}

\begin{proof}
Set  $o =P^{Y} (x)$ and $F=P^Y (F_x)$.
From the equality of slopes we deduce  
$F\times \{ \bar y \} =F_p$ with $p=(o,\bar y) \in Y\times \bar Y =X$. 
This implies that $T_p =T_x$ for each point $p \in T_x$.

For all points $q,p\in T$ the fiber $\bar Y_q$ intersects  $F_p$ in a unique
point.  Moreover, for this intersection  point
$\bar p$ we have $d^2(P^Z (q), P^Z (p)) =d^2 (P^Z (q),P^Z (\bar p))+
d^2 (P^Z (\bar p),P^Z(p))$, due to \lref{interbase}.

 Now consider $Z$ as the union $Z=\cup _{p\in T} P^Z(F_p)$ 
and define the map $P_{pq} : P^Z(F_p)\to P^Z(F_q)$ by
sending $P^Z(\bar p)$ to $P^Z (\bar q)$, where
$\bar q$ is the unique intersection of
$\bar Y _{\bar p} $ with $F_q$. 

The uniqueness of the intersection shows that
$P_{pq} \circ P_{qp} =Id$ and that $P_{pq} \circ P_{rp} =P_{rq}$
for all $r,q,p \in T$. From the above equality we deduce
that for all $z\in P^Z (F_p)$ and all $\bar z \in P^{Z} (F_q)$ we have
$d^2 (z,\bar z ) =d^2 (z, P_{pq} (z)) + d^2 (P_{pq} (z), \bar z)$.

 Hence we can  apply the considerations in Subsection  \ref{recogn} 
 and get the desired conclusion.
\end{proof}

 We believe that the image $P^Z(T)\subset Z$ in the last lemma
always coincides with $Z$. We will prove it below under the
assumption of the finiteness of the affine rank.  First we are going to reduce it 
to the affine case, which then will be finished in the next section.

 We say that $X$
has the property $O$ if for all decompositions 
$X=Y\times \bar Y =Z\times \bar Z$ and each point $x\in X$ the projection
$P^Z :T \to Z$ is surjective, where $T=T_x$ is defined as above.

\begin{lem} \label{propO}
If each finite dimensional affine metric space has property $O$ then
so does each geodesic metric space of finite affine rank. 
\end{lem}

\begin{proof}
Let $X$ be decomposed as $Y\times \bar Y=X=Z\times \bar Z$.
We identify $Z$ with $Z_x$ and use the notations from above.
  Choose a point $z\in Z_x$.
Consider a geodesic $\gamma$ from $x$ to $z$. Let $C$ be a maximal
affine subset of $X$ that contains $\gamma$. Due to \cref{maxexist} such a subset 
exists and due to \cref{maxrect} it is rectangular.
Therefore, $C=A\times \bar A =B\times \bar B$, where $A= P^Y(C), \bar A=
P^{\bar Y} (C) , B = P^Z( C)$ and $\bar B =P ^{\bar Z} (C)$.
 By our assumption $C$ has the property $O$.
 Therefore, for
$\tilde T =P^A (A_x \cap B_x) \times \bar A$ we get
$z \in P^{B_x} (\tilde T)$. Since  
$\tilde T \subset T $ we deduce that $z$ is contained in $P^{Z_x}(T)$.
Since $z$ was chosen arbitrary, we deduce $P^{Z_x} (T) =Z_x$.
\end{proof}

\section{Linear algebra} \label{sec-lin}

\subsection{Euclidean case}
Let $V$ be a finite dimensional vector space. We will write convex
for linearly convex below. We  assume that all convex subsets that appear
 below contain the origin $0$.

 Let $0\in C$ be a convex subset of $V$. We denote by $H(C)$ the linear
  hull of $C$. By $L(C)$ we denote the largest linear subspace that is 
contained in $C$. Observe that $L(C)$ is well defined, actually,  $L(C)$ is
the union of all linear lines in $C$. In particular, $L(C)=0$ if and
only if $C$ does not contain a line.

  For convex subsets $C_1,C_2 \subset V$ their sum is defined as
$C_1+C_2 = \{ c_1 +c_2 | c_i \in C_i \}$. We say that the sum is direct
and write $C_1 \oplus C_2$ if $H(C_1) \cap H(C_2) = 0$.
We say that $C$ is indecomposable if for each decomposition
$C=C_1 \oplus C_2$ one of the summands is $0$.

 We recall the result of Gruber (\cite{gruber}):
\begin{thm}
Let $C$ be a convex subset of $V$ that does not contain  lines. Then $C$
has a unique decomposition $C=C_1 \oplus C_2 \oplus ....\oplus C_k$
with indecomposable $C_j$.
\end{thm}

 Let now $V$ be a Euclidean vector space and $C\subset V$ be a convex subset.
Then  $C$ has a unique decomposition $C=L(C) \oplus C_1 \oplus .... \oplus C_k$,
where all $C_i$ are indecomposable and do not contain lines and where
all $C_i$ are orthogonal to $L(C)$ (compare \cite{gruber}).

 For a decomposition $C=C_1 \oplus C_2$ we have well defined projections
$P^{C_i} :C\to C_i$ defined by $P^{C_i} (c_1+c_2)=c_i$,
for all $c_1 \in C_1, c_2 \in C_2$. We call the decomposition
$C=C_1 \oplus  C_2$ orthogonal if $H(C_1)$ and $H(C_2)$ are orthogonal.
Observe, that the decomposition $C=C_1 \oplus C_2$ is orthogonal if and only
if the natural map $I:C_1 \times C_2 \to C$ given by $I(c_1,c_2) =c_1+c_2$
is an isometry.

\begin{lem} \label{Eucl}
Let $V$ be a Euclidean space and $C$ be a convex subset of $V$. Let
$C=A\oplus \bar A =B \oplus \bar B$ be two orthogonal decompositions.
Then $P^{B} ((A\cap B) + \bar A ) =B$.
\end{lem}

\begin{proof}
Decompose $A$ as $A=L(A)\oplus A_0$, where $A_0$ is orthogonal to $L(A)$.
In the same way decompose $\bar A =L(\bar A) \oplus \bar A_0$,
$B=L(B) \oplus B_0$ and $\bar B=L(\bar B) + \bar B_0$.
We have $L(C)=L(A)\oplus L(\bar A) =L(B) \oplus L(\bar B)$.
Moreover, by the orthogonality assumption 
$A_0 \oplus \bar A_0$ is the orthogonal complement of $L(C)$ in $C$.
Since the same is true for $B_0 \oplus \bar B_0$ we have
$C_0=A_0 \oplus \bar A_0 =B_0 \oplus \bar B_0$.
From the uniqueness of the decomposition of $C_0$ into a direct sum
of indecomposable parts, we deduce $(A_0 \cap B_0) + (\bar A_0 \cap B_0) =B_0$.
In particular, $P^ B ((A\cap B) + \bar A)$ contains $B_0$.

 On the other hand, the equality  $L(B)= L(A) \cap L(B) + P^{L(B)} (L(\bar A))$
is an exercise in linear algebra.
(We may assume
$C=L(C)$.   Considering the quotient space 
$C/((A\cap B) + (\bar A \cap \bar B))$ one is reduced to
the case $A\cap B=\bar A\cap \bar B = \{ 0\}$.  In  this case we must
have $\dim (A)= \dim (\bar B) $ and $\dim (B)=\dim (\bar A)$.
Since the kernel of the projection 
$P^B :\bar A \to B$ is $\bar A \cap \bar B =\{ 0\}$, we deduce that
$P^ B:\bar A \to B$ is injective and therefore surjective as well.)

 Combining the both equalities we arrive at
 $P^{B} ((A\cap B) + \bar A ) =B$. 
\end{proof}

 In the notations used in the last lemma, we have seen 
$B_0 = B_0 \cap A_0 + B_0 \cap {\bar A_0}$. From this we deduce:
\begin{lem} \label{euclhelp}
Let $C=A\oplus \bar A =B\oplus \bar B$ be two orthogonal decompositions
of a convex subset $C\subset V =\R ^n$.  If $B\cap A =B\cap \bar A = \{ 0\}$,
then $B=L(B)$, i.e. $B$ is a linear space in this case.
\end{lem}

\subsection{Banach spaces} Let $V$ be a finite dimensional real vector space.
An ellipsoid in $V$ is the image of the unit ball in some Euclidean space
$\R ^N$  under a linear map $A:\R ^N\to V$.

 Recall that the map that assigns to a norm on $V$  the unit ball $K$
of the norm is a one-to-one correspondence  between norms on $V$ and centrally
symmetric convex subsets of $V$ with a non-empty interior. Recall further,
that the norm stems from a scalar product if and only if the unit ball
$K$ of this norm is an ellipsoid.

 For each norm $|| \cdot ||$ on $V$, there is a unique ellipsoid $E$
of maximal volume that is contained in the unit ball $K$ of $V$
(see for instance \cite{am} or \cite{Th}). This ellipsoid (called the L\"owner ellipsoid)
 defines a scalar product on $V$.  
Hence we have defined an assignment $(V,K) \to (V,E)$ that assigns 
to a norm on $V$ a scalar product on $V$. 
We will denote this Euclidean space arising from a normed space $V$ by $V^e$. 

The following easy observation is probably well known. Since we could not
find a precise reference,we include a short proof:
\begin{lem} \label{ellips}
If a finite dimensional normed vector space $V$ is a direct product
$V=V_1 \times V_2$ of its subspaces $V_1,V_2$ then we have 
$V^e =V_1 ^e \times V_2 ^e$, i.e. $V_1$ and $V_2$ are orthogonal in 
the Euclidean space $V^e$.
\end{lem}

\begin{proof}
Let $K$ be the unit ball of $V$ and let $E$ be the ellipsoid in $K$ of maximal
volume.  Denote by $E_i$ the projection of $E$ to $V_i$ with respect to
the decomposition $V=V_1 \oplus V_2$.   Set 
$\tilde E:= ( E_1 \times E_2)  \cap K$. By construction we have 
$E\subset \tilde E \subset K$.  On the other hand, $E_1$ and $E_2$ are 
ellipsoids  in $V_i$  and $x =(x_1,x_2) \in E_1 \times E_2$ is in $K$
if and only if $||x_1|| ^2 + ||x_2 ||^2 \leq 1$, since $V$ is the direct 
product of $V_1$ and $V_2$.   This shows that $\tilde E$ is the unit ball
of the scalar product on $V$, such that   $V_1$ and $V_2$ are orthogonal with 
respect to this product and such that $E_i$ is the intersection of the
unit ball of this scalar product with $V_i$. 

 By maximality of $E$ we have $E=\tilde E$ and we are done.
\end{proof}

\subsection{General case}
Let now $V$ be a finite dimensional normed vector space.
Let $C$ be a (linearly) convex subset of $V$ with $0\in C$.
Assume that $C$ is decomposed as a direct product $C=C_1 \times C_2$.
Identify $C_i$ with the $C_i$-fiber through $0$. Then $C_i$
is again linearly convex and the assumption $C=C_1 \times C_2$
 is equivalent to the statement that
 $C$ is a direct sum $C=C_1 \oplus C_2$ and that
for all $c_i \in C_i$ we have $||c_1+c_2|| =\sqrt {||c_1||^2 +||c_2||^2 }$.

We claim that $H(C)= H(C_1) \times H(C_2)$. To see this, one can
add a translation and assume that $0$ is an inner point of $C$ in $H(C)$.
In this case, we deduce $||c_1+c_2|| =\sqrt {||c_1||^2 +||c_2||^2 }$
for all $c_i $ in a small neighborhood of $0$ in $C_i$. 
By homogeneity of the norm we deduce 
$||v_1+v_2|| =\sqrt {||v_1||^2 +||v_2||^2 }$ for all $v_i \in H(C_i)$.

Without loss of generality we will assume $V=H(C)$.
 From \lref{ellips} we deduce that if one replaces $H(C)=V $ by $V^e$, i.e. 
if one equips
$V$ with the canonical scalar product defined in the previous subsection
and considers $C$ with the new metric coming from $V^e$, then we still
have $C=C_1 \times C_2$. Observe finally, that
the projections $P^{C_i} :C\to C_i$ remain unchanged under this procedure
since they are defined in purely linear terms. 

 These observations together with \lref{Eucl} show that each finite dimensional
affine metric space has the property $O$, as defined in Section \ref{sec-inter}.
From \lref{lem-proj} and \lref{propO} we deduce:

\begin{cor}  \label{maininter}
Let $X$ be a geodesic metric space of finite affine rank. 
If $Y\times \bar Y =X =Z\times \bar Z$ are two decompositions of $X$, then for
each point $x\in X$ the intersection $Y_x \cap Z_x$ is a direct factor of 
$Z_x$. 
\end{cor}

This reduction of the general affine case to the Euclidean one
and \lref{euclhelp} imply the following:
\begin{cor} \label{unbound}
Let $C$ be an affine metric space with two decompositions
$C=A\times \bar A =B\times \bar B$. If $B\cap A= B\cap \bar A =0$
then $B =L(B)$, i.e. $B$ is a Banach space.
\end{cor}

\section{Reduction}
\label{sec-reduce}
Now we are in position to reduce  \tref{main} to \cref{strike}.

 Observe that each metric space has non-negative affine rank and
each space that contains at least one non-constant geodesic has rank $\geq 1$.
Hence the affine rank of each geodesic metric space with at least two points
is at least $1$.
Due to \cref{additive}, each geodesic
metric space of finite affine rank $m$ can have a decomposition in
at most $m$ non-trivial factors. Therefore each such space has a 
decomposition $X=Y_1 \times Y_2 \times ... \times Y_l$ with irreducible factors
$Y_i$. Rearranging the factors and taking together all factors that are 
isometric   to $\R$, we get at least one decomposition as required in 
\tref{main}.

 It remains to prove the uniqueness. We proceed by induction 
on the affine rank $m$. The case $m=1$ is clear, since in this case 
the space $X$ is irreducible. We assume that the uniqueness holds true
for all $m'\leq m$.  Take a space $X$ of affine rank $m$
and consider a decomposition $X=Y_0 \times Y_1 \times ... \times Y_n$ such
that $Y_0$ is a Euclidean space and such that the $Y_i$ are 
irreducible spaces of positive rank
that are non-isometric to $ \R$, for all $i\geq 1$.
 Let $X=Z_0 \times Z_1 \times ... \times Z_k$ be another such decomposition.  

 We may assume that $X$ is not a Euclidean space and not irreducible,
since the result is clear in these cases.

Fix a point $x\in X$.
 First, assume   $(Z_j)_x =(Y_i)_x$  for some 
$i,j\geq 1 $. Renumerating we may assume  $i=j=1$.  Then
$(Y_0 \times Y_2 \times Y_3 \times  ... \times Y_k) _x =
(Z_0 \times Z_2\times Z_3 \times ... \times Z_l )_x \subset X$, 
since this is the subset
of all points $\bar x$ of $X$  with 
$d(\bar x , (Z_1) _x) =d(\bar x, (Y_1) _x )= d(\bar x,x)$.

 Applying our inductive assumption to the space  
$(Y_0 \times Y_2 \times ...\times Y_k)_x$ we deduce that $k=l$  and that after
a renumeration we have $(Y_i)_x =(Z_i)_x$  for all $i$. From the equality
of slopes (Subsection \ref {equalslope}) we deduce that 
$(Y_i) _{\bar x}  =(Z_i) _{\bar x}$ for all $0\leq i \leq k$ and all 
$\bar x\in X$. 

 From \cref{maininter} and the irreducibility of the $Y_i$ for $i\geq 1$
we know that for all $0\leq i,j \leq k$ we 
either have $(Y_i)_x \cap (Z_j) _x =\{ x \}$ or $(Y_i)_x = (Z_j) _x$.

 Assume now that  the intersection 
$G _x=(Y_0) _x \cap (Z_0) _x$ contains more than one point. Then $G_x$
is a Euclidean space and a direct factor of $(Y_0)_x$ and of $(Z_0)_x$,
due to \cref{maininter}. Write $Y_0$ as $G\times \tilde Y_0$ and
$Z_0$ as $G\times \tilde Z_0$.  Then we have 
$(\tilde Y_0 \times Y_1 \times Y_2 \times ...\times Y_n)_x =
(\tilde Z_0 \times Z_1 \times   ... \times Z_k)_x$.  By our inductive
assumption, we deduce $(Y_i)_x =(Z_j)_x$ for some $i,j \geq 1$.
The argument above shows that the decompositions
$X=Y_0 \times ... \times Y_n$ and $X=Z_0\times ... \times Z_k$ coincide
up to a reindexing.

 Taking both observations together we may assume
that  $(Y_i)_x \cap (Z_j)_x =\{ x \}$ for all $0\leq i\leq n, 0\leq j \leq k$. 

Under this assumptions we set $Y=Y_n$,
$\bar Y=Y_0\times Y_1 \times ... \times Y_{n-1}$, $Z=Z_k$ and
$\bar Z =Z_0 \times Z_1 \times ... \times Z_{k-1}$ and consider
the decompositions $Y\times \bar Y=X= Z\times \bar Z$.

 By assumption we have $Y_x \cap Z_x =\{ x\} $.  If $F_x =\bar Y _x \cap
\bar Z_x$ has more than one point, then $F_x$ is a non-trivial
factor of $\bar Y_x$ and of $\bar Z _x$, due to \cref{maininter}. Take an
irreducible factor $G$ of $F_x$. Applying the inductive assumption
to $\bar Z$ and to $\bar Y$, we deduce that either $G_x =(Y_i)_x =(Z_j)_x$
for some $i,j \geq 1$, or that $G$ is a Euclidean line and that
$G_x \subset (Y_0) _x \cap (Z_0) _x$.  Both conclusions contradict
our assumption.  Therefore $\bar Z _x \cap \bar Y_x =\{ x\}$.

In the same way we deduce $Y_x \cap \bar Z_x = \bar Y _x \cap Z_x= \{ x\}$.
Therefore the decompositions 
$Y\times \bar Y=X= Z\times \bar Z$ satisfy the assumptions of \cref{strike}.
From \cref{strike} (that we will prove in the next section) we deduce
that $X$ is a Euclidean space which is in contradiction with our assumption.


\section{Euclidean rigidity}
\label{sec-ban}
We are going to prove \cref{strike} in this section.  Thus let
$Y\times \bar Y=X =Z\times \bar Z$ be two decompositions of a geodesic
space $X$ of finite affine rank such that 
$Y_x \cap Z_x =Y_x \cap \bar Z_x =\bar Y_x \cap Z_x = 
\bar Y_x \cap \bar Z_x =\{ x \}$.

\subsection{Reduction to Banach spaces}
We claim that $X$ is a Banach space. Indeed, let $C$ be a maximal
affine subset of $X$ that contains the point $x$. Due to \cref{maxrect},
the subset $C$ is rectangular. Hence $A\times \bar A  =C= B\times \bar B$,
where $A=P^Y(C)$, $\bar A =P^{\bar Y} (C)$, $B=P^Z(C)$ and $\bar B=P^{\bar Z} (C)$.
Since $A_x \cap B_x =A_x \cap \bar B_x = \bar A_x \cap \bar B_x  =
\bar A_ x\cap B_x =\{ x\}$, we deduce from \cref{unbound} 
that $B$ and $\bar B$, and therefore also $C$, are Banach spaces. 

  We may assume that $x$ is the origin of the Banach space $C$ and 
identify $A,\bar A, B, \bar B$ with their fibers through $x=0$.
From the transversality of the decompositions we deduce that 
$\dim (A)= \dim (\bar A) =\dim (B) = \dim (\bar B)$. Since the projection
of $A$ to $B$ is injective it is also surjective, i.e.
$P ^Z( A_x )= B_x$. Similarly $P^{\bar Y} (B_x ) =\bar A_x$.
Therefore, each rectangular subset $C_0$ of $X$ that contains 
$A_x$ must contain $C$.

 We claim $\bar A_x = \bar Y_x$. Take an arbitrary point $z\in \bar Y_x$. 
Connect $x$ and $z$ by a 
geodesic $\gamma$. Then $\tilde C=A\times P^{\bar Y} (\gamma ) \subset X$ 
is an affine space.
Take a maximal affine subset $C_0$ of $X$ that
contains $\tilde C$. Then $C_0$ is rectangular (by \cref{maxrect})
and contains $A_x$. Therefore $C_0$ contains $C$ and by maximality of 
$C$ we have $C_0 =C$. Therefore, $z\in C \cap \bar Y_x  =\bar A_x$.
In the same way we see $A_x =Y_x$. Hence $C=X$.

\subsection{Strategy}
Thus we may assume that $X$ is a finite dimensional Banach space.
It remains to prove the following claim:
\begin{claim} \label{claim-banach}
Let $C$ be a Banach space of finite dimension. If $C$ has two decompositions
as a direct product $A\times \bar A= C= B\times \bar B$
such  that $A\cap B =A \cap \bar B =\bar A \cap B =\bar A\cap \bar B = \{ 0 \}$
then $C$ is a Euclidean space.
\end{claim}

 We are going to prove the claim in the following manner. If $B$ and $A$
enclose a unique angle, i.e. if each line of
$B$ has the same slope with respect to $A$, then the relation between
the both decompositions induce an easy condition on the norm, that 
turns out to be equivalent to the parallelogram equality 
(Subsection \ref{subsec-comp-projections}). For the general case, we prove in
the remaining part of this section that there
is a rectangular subspace $L$ of $C$ that is as non-Euclidean as $C$
and such that its decompositions ``enclose a unique angle''. This subspace $L$
is constructed as the projection of some extremal subsets
in the product space $C^2$, the extremality being described as the maximal
possible violation of the parallelogram equality.

 Before we embark on the proof, note that
a subspace $V$ of a Banach space $C$ is rectangular with
respect to a decomposition $C=C_1 \times C_2$ if and only if
$V= V\cap C_1 + V\cap C_2$.  This implies that $V+W$ is
rectangular with respect to the decomposition $C=C_1 \times C_2$
if $V$ and $W$ are rectangular.


\subsection{Projections and their compositions}
\label{subsec-comp-projections}
Let $C$ be decomposed as above. Consider the bijective linear
projections $P^A:B\to A$ and $P^B:A\to B$, and let $Q:=P^A \circ P^B :A\to A$
be their composition.

\begin{lem} \label{lemma-eigenvector}
Under the above assumptions $Q$ has an eigenvector.
\end{lem}
 
\begin{proof}
Let $x\in A$ be arbitrary. The vectors $P^B(x)$ and 
$P^{\bar B} (x)$ span a two-dimensional Euclidean subspace $F$
of $C$ that contains
$x$. If we denote by $\alpha $ the angle between $x$ and $P^B(x)$ then by
definition $||P^B(x)|| =\cos(\alpha ) ||x||$. On the other hand, the projection
$\bar x$ of $P^B(x)$ to the line spanned by $x$ in $F$ satisfies 
$||\bar x|| = ||P^B(x) || \cos (\alpha )$. This implies
$\frac {||P^B(x)||} {||x||} \leq \frac {||P^A (P^B (x))||} {||P^B(x)||}$ 
and equality holds if and only if $x$ is an eigenvector of $Q$.

Similarly, we have $\frac {||P^A (\bar x)||} {||\bar x||}  \leq
\frac {||P^B (P^A (\bar x))||} {||P^A(\bar x)||}$ for all $\bar x \in B$.
Thus, taking a point $x$ in the unit sphere of $A$ such that
$\frac {||P^B (x)||} {||x||}$ is maximal, we deduce that $x$
must be an eigenvector of $Q$.
\end{proof}

 Note that if some $a\in A$ satisfies $Q(a)= \lambda a$, then
for $b= P^B (a)$ we have $(P^B \circ P^A) (b) =\lambda b$.
Moreover, we have $P^A \circ P^{\bar B}  (a)= (1-\lambda ) a$.

  This shows that if $a\in A$ is an eigenvector of $Q$,
then $P^B(a)$ and $P^{\bar B} (a)$ generate a two-dimensional 
Euclidean subspace $F$ of $C$ that is rectangular with
respect to both decompositions of $C$. 
Making the computations in the Euclidean space $F$, 
one concludes that $0< \lambda <1$ and that  
 $||P^B(a)|| = \sqrt \lambda ||a||$ holds.

\begin{lem} \label{unique}
Under the above assumptions let in addition $Q$ be a multiple of the identity:
$Q=\lambda \Id _A$. Then $C$ is Euclidean.
\end{lem}

\begin{proof}
We have seen that in this case $P^{\bar A} \circ P^{B} =(1-\lambda ) 
\Id_{\bar A}$. 
Moreover, $P^A:C\to A$ satisfies
$||P^B (a)||=  \sqrt \lambda ||a||$ for all $a\in A$.
In the same way, $||P^{\bar A} (b)|| =\sqrt {1-\lambda } ||b||$ for all
$b\in B$.  

Consider the map $I:A\to \bar A$ given by 
$I(a) = \frac 1 {\sqrt {\lambda (1-\lambda )}}   (P^{\bar A} \circ P^B) (a)$.
Then $I$ is an isometry. Identify $A$ and $\bar A$ via this isometry.
Then the space $C$ is identified  with $A^2=A\times A$ and 
the subspaces $B$ and $\bar B$ of $C$ are given as
$B= \{(x,sx)|x\in A \}$ 
 and $\bar B = \{(-sx,x)|x\in A\}$, where 
with $s=\sqrt{\frac{1-\lambda}{\lambda}}$.

The assumption $A^2 =C= B\times \bar B$ now reads as
$||(x,sx)||^2 + ||(-sy,y)||^2 = ||(x-sy,sx+y)||^2$, for all $x,y\in A$.
Thus we have $(1+s^2)||x||^2+ (1+s^2)||y||^2)=||x-sy|| ^2 +||sx+y||^2$
for all $x,y \in A$. 

  For $s=1$ this is just the usual parallelogram equality which implies
that $A$ is Euclidean. For general $s\neq 0$ it is shown in
\cite{carl} (compare also \cite{am}, 1.16) that this equality also implies
that $A$  is Euclidean. Since $C=A\times A$, the conclusion follows.
\end{proof}


\subsection{Projections and squares}
\label{subsec-squares}
Let a Banach space $C$ be decomposed as $C=A\times \bar A$.
Consider the Banach space $C^2 =C\times C$ and its decomposition
$C^2 =A^2 \times \bar A ^2$.  Let $L$ be a subspace of $C^2$ that
is rectangular with respect to this decomposition. Then for each
point $(v,w) \in L \subset C^2 $ 
there are unique  points $(v_1,w_1) \in L \cap A^2$
and $(v_2,w_2) \in L\cap \bar A ^2$ such that $(v,w)=(v_1,w_1)+(v_2,w_2)=
(v_1+v_2,w_1+w_2)$.  Denote by $V$ and  $W$  the projection
of $L \subset C^2$ onto the first and the second $C$-factor, respectively
(i.e. $V$ is the set of all $v\in C$, such that for some $w\in C$, the point
$(v,w) \in C^2$ is contained in $L$).
Then $V,W\subset C$ satisfy $V=V\cap A +V\cap \bar A$ and
$W=W\cap A +W\cap \bar A$. Therefore, $V,W$ and $V+W$ are rectangular
subspaces  with respect to $C=A\times \bar A$.

Assume now that $C$ has two transversal
 decompositions $A\times \bar A =C =B\times \bar B$ 
as above.  Consider the corresponding
decompositions $A^2\times \bar A^2= C^2 =B^2 \times \bar B ^2$ and let
$L\subset C^2$ be a subspace that is rectangular  with respect to  both
decompositions.  Let $V$ ($W$, resp.) be the projection of $L$ onto the first
(the second, resp.) $C$-factor of $C^2 =C\times C$.

Consider the orthogonal projections $P^{A^2} :B^2 \to A^2$
and $P^{B^2} :A^2 \to B^2$ and the composition
 $\tilde Q = P^{A^2} \circ P^{B^2}$.

In the same way define $Q_V:V\cap A \to V\cap A$  and $Q_W :W\cap A\to
W\cap A$. From the definition we conclude that 
the equality $\tilde Q=\lambda Id _{L\cap A^2}$ for some $\lambda \in \R$ 
implies that $Q_V =\lambda Id_{V\cap A}$ and $Q_W =\lambda Id_{W\cap A}$. 

In this case we see that $(V+W)\cap A$ is contained in the 
$\lambda $-eigenspace of $P^A\circ P^B:A\to A$. From  \lref{unique}
we deduce  that $V+W$ is a Euclidean space.


\subsection{Extremal points}
Recall that a Banach space $C$ satisfies the parallelogram
inequality $||x+y||^2 +||x-y||^2 \leq 2(||x||^2 +||y||^2)$ for all $x,y\in C$
if and only if $C$ is Euclidean. In this case the
above inequality is in fact an equality for all pairs $(x,y)\in C ^2$.

We may reformulate this condition as follows. Consider the Banach space
$C^2$ with the product norm and the linear map 
$D:C^2 \to C^2$ defined by 
$D((x,y))=\frac 1 {\sqrt 2} (x+y,x-y)$ for all $x,y\in C$.
Then $C$ is a Euclidean space if and only if the bijection $D$
has norm $1$. In this case $D$ is in fact an isometry.

We define  $M(C)$ to be the norm of the linear map $D$ and we denote by
$E(C)$ the set of all  $v\in C^2$ for which
$||D(v)||=M(C) ||v||$ holds.  By definition $\lambda v $ is contained
in $E(C)$ for all $v\in E(C)$ and all $\lambda \in \R$.

Let now $C$ be decomposed as $C=A\times \bar A$ and consider the induced
decomposition $C^2 =A^2 \times \bar A ^2$. 
The subsets $A^2$ and $\bar A ^2$ are invariant under the map $D$
defined above. Therefore  for arbitrary $v\in A^2$ and
$w\in \bar{A}^2$, we get $||D(v+w)||^2 = ||D(v) +D(w)||^2 =
||D(v)||^2 +||D(w)||^2 \leq M(C) ||v||^2 + M(C) ||w|| ^2= M(C) ||v+w||^2$.
Moreover, equality holds if and only if $v$ and $w$ are contained in $ E(C)$.

This shows that the subset $E(C)$ of $C^2$ is rectangular with 
respect to the decomposition $C^2 =A^2 \times \bar A ^2$.


\subsection{Extremal subspaces}
Now we are in position to finish the proof of Claim \ref{claim-banach}.
Thus let $A\times \bar A= C=B\times \bar B$ be two transversal
decompositions.
 Let  $M(C)\in \R ^+$ and 
$E(C)\subset C^2$ be defined as in the last subsection.

 Choose a largest linear subspace $L_0$ of $C^2$ that is contained in
$E(C)$.  By definition $E(C)$  contains at least one line, hence
$\dim (L_0) \neq 0$.  Since $E(C)$ is rectangular with respect 
to the decompositions $C^2 =A^2 \times \bar A ^2$ and
$C^2 =B\times \bar B^2$, the largest subspace $L_0$ of $E(C)$
is rectangular with respect to these decompositions as well.

 Denote by $\tilde Q: L_0 \cap A^2 \to L_0 \cap A^2 $ the composition
of  projections 
$\tilde Q= P^{A^2} \circ P ^{B^2}$.  Due to \lref{lemma-eigenvector} there is
an eigenvector $v$ of $\tilde Q$. Then the two-dimensional
subspace $L$ of $L_0$ that is generated by 
$P^{B^2} (v)$ and $P^{\bar B^2} (v)$ is rectangular with respect
to the decompositions $A^2\times \bar A ^2 =C^2 = B^2 \times \bar B^2$.

 Since $L\cap A^2$ is one-dimensional, the restriction of
$\tilde Q$ to $L\cap A^2$ is a multiple of the identity.
As in Subsection \ref{subsec-squares}, denote by $V$ and $W$ the projections
of $L\subset C\times C$ onto the first and the second $C$-factor, respectively.
Finally, set $\tilde C= V+W$.  In  Subsection \ref{subsec-squares}
we have seen that $\tilde C$ is a Euclidean space.

 On the other hand, $L\neq \{ 0\}$, hence either $V$ or $W$ are not $\{ 0\}$.
Without loss of generality assume $V\neq \{ 0\}$ and choose some $v\in V$,
$v\neq 0$. By construction there is some $w\in W$ with
$(v,w) \in L \subset E(C)$. Hence $||D((v,w))||=E(C)||(v,w)||$ and,
since $(v,w) \in \tilde C ^2 \subset C^2$ we deduce that
$E(C)=E(\tilde C)$.   But $\tilde C$ is Euclidean, hence $E(\tilde C) =1$.
Thus $E(C)=1$ and $C$ is a Euclidean space.

\bibliographystyle{alpha}
\bibliography{prod}

\begin{thebibliography}{{Mos}92}

\bibitem[Ami86]{am}
D.~Amir.
\newblock {\em Characterizations of Inner Product Spaces}.
\newblock Birkh\"auser Verlag, Basel, 1986.

\bibitem[Bau75]{baum}
G.~Baumslag.
\newblock Direct decompositions of finitely generated torsion-free nilpotent
  groups.
\newblock {\em Math. Z.}, 145:1--10, 1975.

\bibitem[Car62]{carl}
S.~Carlsson.
\newblock Orthogonality in normed linear spaces.
\newblock {\em Ark. Mat.}, 4:297--318, 1962.

\bibitem[dR52]{rham}
G.~de~Rham.
\newblock Sur la reductibilit\'e {d'un} espace de {Riemann}.
\newblock {\em Comment. Math. Helv.}, 26:328--344, 1952.

\bibitem[EH98]{eh}
J.~Eschenburg and E.~Heintze.
\newblock Unique decomposition of {Riemannian} manifolds.
\newblock {\em Proc. Amer. Math. Soc}, 126:3075--3078, 1998.

\bibitem[Fou71]{four}
G.~Fournier.
\newblock On a problem of {S. Ulam}.
\newblock {\em Proc. Amer. Math. Soc}, 29:622, 1971.

\bibitem[FS02]{fs}
T.~Foertsch and V.~Schroeder.
\newblock Minkowski versus {Euclidean} rank for products of metric spaces.
\newblock {\em Adv. Geom.}, 2:123--131, 2002.

\bibitem[Gru70]{gruber}
P.~Gruber.
\newblock Zur {Charakterisierung} konvexer {K\"orper}. \"uber einen {Satz} von
  {Rogers} und {Shephard}.
\newblock {\em Math. Ann.}, 184:79--105, 1970.

\bibitem[Her94]{intui}
I.~Herburt.
\newblock There is no cancellation law for metric products.
\newblock In {\em Intuitive geometry (Szeged, 1991)}. North-Holland, Amsterdam,
  1994.

\bibitem[HL]{hitz}
P.~Hitzelberger and A.~Lytchak.
\newblock Spaces with many affine functions.
\newblock Preprint: Arxiv: \emph{math.MG/0511583} To appear in Proc. Amer.
  Math. Soc.

\bibitem[{Mos}92]{mos}
M.~{Moszy\'nska}.
\newblock On the uniqueness problem for metric products.
\newblock {\em Glas. Mat. Ser. III}, 47:145--158, 1992.

\bibitem[Tho96]{Th}
A.~C. Thompson.
\newblock {\em Minkowski Geometry}.
\newblock Cambridge University Press, Cambridge, 1996.

\end{thebibliography}

\end{document}